\documentclass{article}

\usepackage{amsmath,amsthm,amsfonts}

\usepackage{amssymb}

\usepackage{epsfig}

\usepackage{rotating}

\usepackage{subfigure}

\def\sl R{\mathscr{R}}

\begin{document}

\title{Numerical scaling analysis of the small-scale structure in turbulence }
\author{Panagiotis Stinis and Alexandre J. Chorin\\ 
\\
Department of Mathematics \\
   University of California \\
and \\
Lawrence Berkeley National Laboratory \\
   Berkeley, CA 94720}
\date {}

\maketitle

\begin{abstract}
We show how to use numerical methods within the framework of successive scaling to analyse the 
microstructure of turbulence, in particular to find inertial range exponents and structure functions. The 
methods are first calibrated on the Burgers problem and are then applied to the 3D Euler equations. 
Known properties of low order structure functions appear with a relatively small computational outlay; however, 
more sensitive properties cannot yet be resolved with this approach well enough to settle ongoing controversies.
\end{abstract}

\section{Introduction}\label{intro}
It has often been suggested that scaling and renormalization group ideas should be helpful in the analysis of the 
small scales of turbulence at high Reynolds numbers, and there is a vast literature in the subject but few 
concrete results (see e.g. \cite{smith} and references therein). One of the outstanding difficulties is that scaling ideas are usually 
implemented within a perturbative framework, and in the absence of a small parameter the validity of the 
framework remains doubtful.

In an earlier paper \cite{chorin12}, one of us has attempted to marry scaling with a particular numerical method; 
the results could hardly be called definitive. In the current paper we try again in the context of a spectral method.
For reviews of related methods in turbulence and statistical physics, see \cite{goldenfeld,barenblatt}.
Related ideas have also been presented in \cite{berger,landman,frederico}.

Suppose you can represent on the computer Fourier modes up to the some wave number $N$.  All the equations we work with in the 
present paper contain an energy 
cascade; when the energy reaches wavenumber $N$ aliasing begins, energy is reflected into the longer wavelengths in ways 
that are not justified by the equations, and the approximation becomes invalid. On the other hand, the characteristic time of the 
modes becomes shorter and one could view the amplitudes of larger modes as nearly stationary on the time scales of the smaller modes. 
When the energy reaches mode $N$ one could think of rescaling the computation so that small scale modes are added, large scale modes are 
removed from the computation because they are nearly constant, until a new cutoff $N_1$ is reached, and so on, thus computing 
in a moving window of modes and probing the large wavenumber coefficients of the flow. The problem is how to do this scaling 
and how to justify the results.

We first explore the rescaling method in the case of Burgers equation where the structure of the flow is well 
understood, and we show how to pick the rescaling time and how to rescale the  flow; we then verify that well-known results are 
reproduced. The key difficulty, as in some other numerical renormalization methods \cite{gear}, has to 
do with the rescaled boundary conditions.

We then apply the idea to the Euler equations in 3D; well-known results are rediscovered at low cost, attempts at settling 
current controversies are made, and ideas for future improvement are suggested. The paper represents work in progress; we felt that the 
methodology is promising and worthy of presentation. 
What should be done next is discussed in the final section.

We focus in particular on structure functions, which are averages of the velocity field of the form $S_n(r)=<(u(x+r)-u(x))^n>$, where
$x, x+r$ are points in the fluid $r$ apart, $u$ is the velocity at this points, $n>0$ is a power, and the brackets denote
a (spatial, temporal, or ensemble) average. Some of the key results in turbulence theory relate to these functions.
The Kolmogorov ``K41" theory \cite{kolm} deduced that $S_n(r)=C_nr^{n/3}$, where $C_n$ is a constant that depends on $n$ only.
This would be an exact result if the velocity field were Gaussian; however, the velocity field in
turbulence is not Gaussian (see e.g. \cite{batchelor,bernard}). It is well-known that experiment gives exponents different than the
Kolmogorov values (\cite{monin}). The exponent for $n=2$ in particular has given rise to much controversy. The case $n=3$ is special
because the conclusion $S_3=C_3r$ is ``almost" a theorem (\cite{monin}). In recent years Barenblatt et al. \cite{bcbc,BC04}
have conjectured that the structure funnction exponents may be Reynolds-number dependent.

Note that in some sense the scaling transformations here accomplish the opposite of what is usually attempted
with real-space renormalization methods (as in \cite{kadanoff, goldenfeld, CH}): the focus here narrows to ever smaller
scales rather than expand to ever larger scales. The point is of course
that in either case one uses scaling trasnformations to explore those properties of a system that are scale invariant.

The paper is organized as follows. Section \ref{algo} contains  an explanation of the construction; section \ref{burg} describes a validation 
in the well-understood case of the one-dimensional inviscid Burgers equation. 
Section \ref{eulns} contains results about the Euler equations. A discussion  and ideas about futher work 
follow in Section \ref{conc}.


\section{The scaling algorithm}\label{algo}
We present the algorithm in the case of the Navier-Stokes equations in three space dimensions. The modifications needed for 
the Euler equations, fewer dimensions, and in the case of the Burgers equation are straightforward. Consider therefore the 
3D Navier-Stokes equations with periodic boundary conditions in the box $[0,2\pi]^3$:

\begin{gather}
\label{ns}
u_t+u\cdot \nabla u= - \nabla p + \nu  \Delta u ,\; \nabla \cdot u=0,  \\
u(x,0)=u_0(x), \quad u(x+2\pi e_i)=u(x), \, i=1,2,3, \notag
\end{gather}
where $u(x)=(u_1(x_1,x_2,x_3),u_2(x_1,x_2,x_3),u_3(x_1,x_2,x_3)),$ 
$\nabla= (\frac{\partial}{\partial x_1},\frac{\partial}{\partial x_1},\frac{\partial}{\partial x_1} ),$ 
$\Delta = \sum_{i=1}^3\frac{{\partial}^2}{\partial x^2}$ and $\nu$ is the viscosity. 
Also, $e_1=(1,0,0)$ and similarly for the other two directions.

We proceed in the following sequence of steps:

\begin{enumerate}

\item
Prescribe an initial condition.

\item
Solve the system of ordinary differential equations for the Fourier modes. If a certain threshold (to be specified below), 
is crossed, indicating the onset of aliasing,  stop the calculation.

\item
Transform to real space and locate the point where the vorticity $\omega= \nabla \times u$ 
is largest.

\item
Construct a box of size $\pi^3$ centered at the point where the vorticity 
is largest. This ``reduced" box contains half the  points in 
each spatial direction; equivalently, the solution 
in this reduced box can be described by eighth of the number of Fourier modes needed for describing the solution in the full box (this accounting will be slightly modified after we construct boundary conditions, see below).

\item
Rescale quantities and parameters as needed while stretching the solution in the reduced box to the whole box. Transform the stretched solution 
to Fourier space. Fewer Fourier modes are needed to describe the solution 
than before the reduction and stretching. Fill the rest of the available Fourier modes 
with zeros.

\item
Use the stretched solution as initial condition and repeat steps 2-5. 
Steps 2-5 will be referred to from now on as a cycle.

\end{enumerate}

A sequence of such cycles will do what we promised to do in the introduction: compute in a 
moving window of modes while focusing on a singularity.

We now summarize some of the difficulties that must be addressed in all cases. 
The major one is that the velocity in the reduced box is not periodic.
If one continues it periodically so that the computation proceed one creates discontinuities in every period. 
These discontinuities can dominate the solution and produce harmful artifacts. Similar difficulties
with boundary conditions appear in other implementations of renormalization ideas (see e.g. \cite{gear}).
In addition, in several space dimensions
one has to take as initial condition not the stretched flow in the reduced box but its 
divergence-free part, so that the initial conditions satisfy the equation of continuity. 
In the present paper we minimize the harm by a ``fringe" construction 
(see e.g. \cite{spalart,nordstrom}). We add a ``fringe"  at the end of the reduced domain, 
inside which the solution is forced to vary smoothly, and then impose periodic boundary conditions on solution in the extended domain (reduced domain + fringe). 
In the case of the Euler equations one has to use a three-dimensional fringe region. 
We end up with a box of size $(\frac{3\pi}{2})^3.$ This is our extended reduced box which has to be stretched etc. The 
extended box must always be smaller than the original box. 
To minimize the impact of the fringe method on the computation one adds a forcing
term to the equations in the fringe region; for details, see \cite{spalart,nordstrom}.

In a more extended paper to come we will exhibit various treatments of the
boundary conditions and their effects on the results, and see that the results
do not depend sensitively on what is done at the boundaries. However, this
is clearly not the end of the story. The derivation of an optimal matching of
the flows at the different scales to each other so as to miminize boundary effects
is a problem we plan to investigate further in the future.

The criterion we 
use to decide that aliasing is becoming dangerous and  that one should rescale is based on
   the ratio of the energies in the outermost and the innermost 
shell in Fourier space, i.e., the ratio of the energy in the 
smallest scales and the energy contained in the largest scales. If the ratio becomes larger 
than a prescribed tolerance $\epsilon$, we stop the calculation and  rescale. 
The choice of $\epsilon$ will also be discussed in detail below.
In addition, during each cycle we monitored the evolution of the Taylor scale 
\begin{equation}
\lambda=\biggl(  \frac{5\int_{0}^{2\pi}\int_{0}^{2\pi} \int_{0}^{2\pi}
   |u|^2 dx_1dx_2dx_3} {\int_{0}^{2\pi}\int_{0}^{2\pi} \int_{0}^{2\pi}
   |\nabla u|^2 dx_1dx_2dx_3} \biggr)^{\frac{1}{2}}.
\label{taylor}
\end{equation}
The Taylor scale signifies the length scale above which there is significant amount of of energy. Thus, in the calculations, this 
scale should always remain larger than the smallest resolvable scale. 
During all the cycles, the Taylor scale never became less than twice the smallest resolvable scale. 
This is enough, given that the algorithm focuses on the regions with large gradients, thus 
overestimating the denominator in (\ref{taylor}).

For ease of computation it is necessary to rescale the variables and 
parameters appearing in the equations after each scaling. 
Each new start corresponds to a change of variables, which has to be performed 
in a way that respects the equations of motion
and keeps all of the variables within a range that the computer can handle with ease. Consider the equations (\ref{ns}) and perform the change of variables 
$$x'=\alpha x, \; t'=\beta t, \;  u'=\gamma u , \;  p'=\delta p.$$ The equations become
\begin{equation}
\label{nsres}
\frac{\beta}{\gamma}{u'}_{t'}+\frac{\alpha}{\gamma^2} u' \cdot \nabla' u'= 
- \frac{\alpha}{\delta}\nabla' p' + \nu \frac{\alpha^2}{\gamma} \Delta' u'. 
\end{equation}
The primed variables are the variables after the stretch. If we use a fringe region, 
then $\alpha=\frac{4}{3}$. We chose $\beta=\frac{4}{3}$ so that the velocities have numerical values
within a fixed range; $\gamma=1, \delta=1$ so that the terms in the equation other
than the viscosity term do not acquire new coefficients as the rescalings proceed.

Equations (\ref{nsres}) become
\begin{equation}
\label{nsres2}
{u'}_{t'}+u' \cdot \nabla' u'= 
- \nabla' p' + \frac{4}{3}\nu  \Delta' u' 
\end{equation}
With our choice of rescaling parameters, we end up 
with an equation where the viscosity increases by a factor of $\frac{4}{3}$ at each rescaling (as one may expect from the fact that the 
viscosity acts more strongly on smaller scales).

A time $t_n$ in the calculation after $n$ scalings corresponds
to an elapsed time of $t_n/(\frac{4}{3})$ in the calculation after only $n-1$ scalings.
If $t_i, i=1,\ldots,M$ is the time between the $(i-1)$-the scaling and the $i$-th scaling,
(the zero-th scaling being the absence of scaling at the beginning of the whole process),
the total time elapsed in real time is $\sum_{i=0}^{M} t_i/(\frac{4}{3})^i$. 
Suppose a singularity forms in the flow in a (real) time $T$.
The support of the spectrum will then be unbounded; the window in which one computes
will move to infinity, requiring an infinite number of scalings; the quantity
$\sum_{i=0}^M t_i/(\frac{4}{3})^i$, if it converges, is an estimate of the singularity formation time $T$,
with the obvious caveats involving the unknown convergence properties of our
scheme.

Finally, note that in more than one space dimension it is ambiguous 
what the largest value of the vorticity is; it may be the largest value of the
modulus of the vorticity or the largest component of the vorticity \cite{xconst}. We made
runs in which one, then the other, were used in the determination of the
center of the rescaled box, and found that this made no difference. For definiteness,in the runs 
below we assume that the point where the vorticity is largest is the
point where one of the components is largest. 

\subsection{The inviscid Burgers equation}\label{burg}

We validate the algorithm above by applying it to the inviscid Burgers equation with an initial condition that gives 
rise to a shock wave whose position and time of occurence can be found analytically. Consider the inviscid Burgers 
equation with periodic boundary conditions in $[0,2\pi]$ 
\begin{equation}
\label{burgers}
u_t+u u_x= 0 
\end{equation}
and the initial condition $u(x,0)=\cos(x)$,  which gives rise at time $T=1.$ to a shock wave located at $x=\frac{\pi}{2}$.
We use this fact to pick a scaling criterion and 
calibrate the algorithm. 
We use $N$ Fourier modes 
($\frac{N}{2}$ positive and $\frac{N}{2}$ negative) to resolve the solution; we rescale and restart when 
$(|u_{\frac{N}{2}}|/|u_{1}|)^2 \geq \epsilon.$ 
We want to find the value of $\epsilon$ for which the total time 
$\sum_{i=0}^{\infty}t_i/(\frac{4}{3})^i$ approximates the known value of $T$.
We approximate
the sum $\sum_{i=0}^{\infty}t_i/(\frac{4}{3})^i$ by $\sum_{i=0}^{45}t_i/(\frac{4}{3})^i$; after 45 scalings, 
the time spent in a cycle has shrunk down to about $10^{-8},$ so that this is acceptable.

We present results from a calculation with N=1024; the equations were solved by a
Runge-Kutta-Fehlberg method  (\cite{hair}) with the tolerance per unit step set to $10^{-10}$ and the $3/2$ rule \cite{canuto} was
used for dealiasing. The energy ratio
restarting criterion was $\epsilon= 2.5\times 10^{-7} $, which yields a total time for the 45 
cycles of 1.014, a $1\%$ error.

The structure functions are given by
\begin{equation}
S_n(r)=\frac{1}{2\pi}\int_{0}^{2\pi} (u(x+r)-u(x))^n dx, \; r  \in [0,2\pi]
\label{struct1d}
\end{equation}
for $n=2,3,4,5.$ If the structure functions are invariant under scaling transformations, then 
the structure functions for the different cycles should have the same form. Plotted in log-log 
coordinates in the original scale, the structure functions for all cycles should be 
translates of each other. Equivalently, though each cycle 
operates on a scale which is half that of the previous cycle, common features 
in the structure functions for different cycles should appear  when the structure functions for 
all cycles are plotted for arguments in $[0,2\pi].$ Finally, averaging the structure functions 
for the different cycles should bring out these common features. 
We present results for 
structure functions averaged over the several cycles. 
The variance of the results provides some rough measure of the uncertainty, and is
very small in the present case of a well-resolved one-dimensional calculation. Also, note that the 
structure functions can only be accurate for small distances because the solution differs for large distances 
from cycle to cycle. Thus, averaging over the different cycles can lead to an erroneous estimate of the 
structure functions for large distances. The same appplies for the case of the Euler equations.

\begin{figure}
\centering
\epsfig{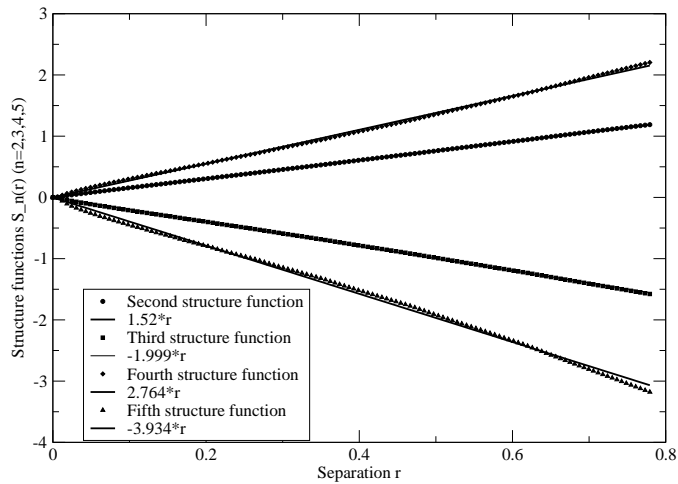}
\caption{Structure functions of orders 2-5 averaged over 45 cycles for inviscid Burgers equation with 1024 Fourier modes}
\label{fig_bur2}
\end{figure}

The structure functions  $S_n(r)$ for Burgers equation are proportional to $r$ for all 
orders $n$ greater than 1, with odd order structure functions negative for small $r.$ 
Figure \ref{fig_bur2} shows the average of the structure functions (\ref{struct1d}) over 45 cycles. 
We do not include the 0th cycle in the averaging of the structure functions because the zeroth step is an 
an equilibration step where the specifics of the initial condition dominate. The averages of the structure functions do exhibit this linear behavior 
and the value of the slope is accurate to three digits. An alternative way to compute the structure functions is by integrating over a 
restricted interval centered at the shock. We performed this computation and the results did not change for small distances (see also comment 
in the previous paragraph).


\section{The Euler equations}\label{eulns}
We now present results for the structure functions of the Euler equations with the Taylor-Green vortex as initial condition:
\begin{gather*}
u_1(x,0)=\sin(x_1)\cos(x_2)\cos(x_3), \\
u_2(x,0)=-\cos(x_1)\sin(x_2)\cos(x_3),\\
u_3(x,0)=0
\end{gather*}
For the Euler equation we cannot afford 
a resolution of $N=1024^3$ without parallelization. We 
plan to report on results  with such resolution in the future. 
Here we present results for a resolution of 
$N=32^3$. The numerical method used 
is the same as in the case of Burgers with the error tolerance set to 
$10^{-10}.$

There are differences here from the 
Burgers case. We assume that the spectrum at high frequencies is dominated
by regions of high vorticity (see e.g. \cite{chorin4}).
This assumption is not universally accepted. We need make no assumption about the structure
of this high-vorticity region, indeed, we hope that a method such as ours can eventually
reveal it. 
In the absence of such assumption, the
value of the cycle restarting criterion $\epsilon$ has to be determined empirically. We used as an initial guess,the value 
$\epsilon=5\times 10^{-4}$ from the 
Burgers case with $N=32$ ; this is the value which in the Burgers case, guarantees 
that a shock wave which develops when the initial condition is $u(x,0)=\cos x,$ is captured at
the right time. We then 
experimented with larger  and smaller values of $\epsilon$ and looked at the form of the structure 
functions. 
Our experiments, with different values of $\epsilon$, show that there is a core range 
of scales where the structure functions exhibit power law behavior. This range can disappear only if the 
calculation is underesolved (large $\epsilon$) or the cycles are restarted too soon for the dynamics to shape the 
structure functions (small $\epsilon$).
The results we present are for the criterion value $10^{-4}$ which 
gives rise to power law behavior for approximately $0 \leq r \leq 1.$ 
The energy ratio criterion measures the ratio of the energy in the outer shell (in Fourier 
space) to the inner shell (a few more comments about monitoring the resolution of the calculation are 
given below).

In addition, in three dimensions, the large gradient structures 
are not spatially localized. Our algorithm centers around the point of
largest vorticity, thus it will inevitably lead to a chopping of structures that protrude from the reduced 
box on which we focus. This decreases the range of scales for which the results can be indepenedent 
of the boundary conditions. In other words, the chopping of the structures affects the large scales 
of the solution. However, we find that it does not afffect the smale scale structure. The structure functions are given by 
$$S_n(r)=\frac{1}{({2\pi})^3}\int_{0}^{2\pi}\int_{0}^{2\pi} \int_{0}^{2\pi}
   (u(x+r)-u(x))^n dx_1dx_2dx_3, \; r \in [0,2\pi]$$
for $n=2,3,4.$

Each time we rescale and restart a calculation we need boundary
conditions at the edge of the then new box as in the Burgers case. As before, we decrease
the damage by the fringe construction. 
For the problem to make sense the new initial conditions have to be divergnce-free
with the new periodic boundary conditions, so the new  initial conditions 
with the periodic boundary conditions must be projected \cite{chorin0} on the space of divergence
free vector fields. We do not expect the projection step to harm the spectrum
or the dynamics, because the projection is equivalent to the addition of the
harmonic vector field and leaves the vorticity invariant.

\begin{figure}
\centering
\subfigure[]{\epsfig{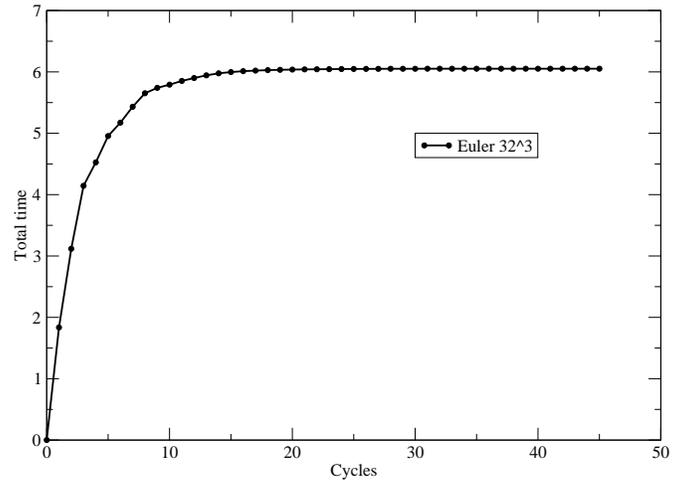}}
\qquad
\subfigure[]{\epsfig{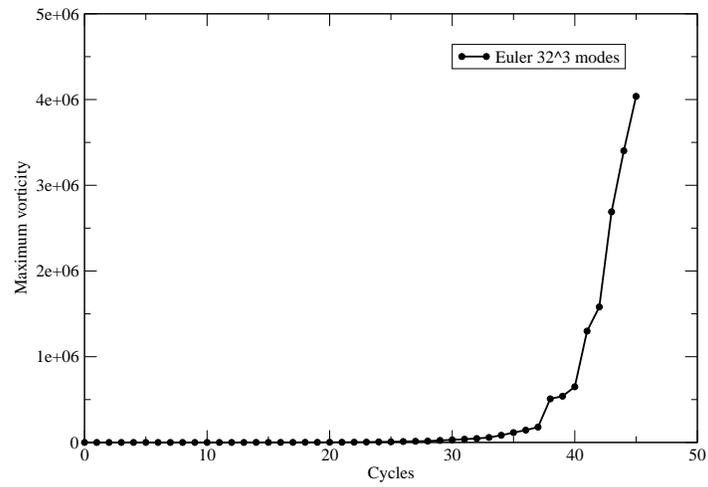}}
\caption{Total time and maximum vorticity evolution with cycles for the Euler equations with $32^3$ modes. }
\label{fig_eul2}
\end{figure}

We begin with results regarding  the possibility of finite time blow-up. Figure \ref{fig_eul2}a shows the evolution of the 
total time as a function of cycle with $\epsilon=10^{-4}.$ Figure \ref{fig_eul2}b shows the evolution of 
the maximum vorticity for the Euler equations as a function of the number of cycles. We 
stopped our calculations after 45 cycles. 
The time spent in the 45th cycle was of $O(10^{-5}).$ It is evident that the vorticity grows in a manner 
consistent with a finite-time blow up around $T=6.05$. When 
the 45th cycle is over, the maximum vorticity has grown by a factor of $10^{6}$. 
Clearly one cannot use this fact to conclude that the vorticity indeed blows up, because
of the uncertainty over the effect of the boundary conditions;  note that
the approximate blow-up time for the Green-Taylor initial
data based on Taylor expansion of the solution \cite{morf} gave (an equally unreliable) blow-up time of T=5.2. 
With $\epsilon=10^{-5}$ the estimate of T decreased to 5.12, while with $\epsilon=10^{-3}$ it increased to 9.04.

\begin{figure}
\centering
\subfigure[]{\epsfig{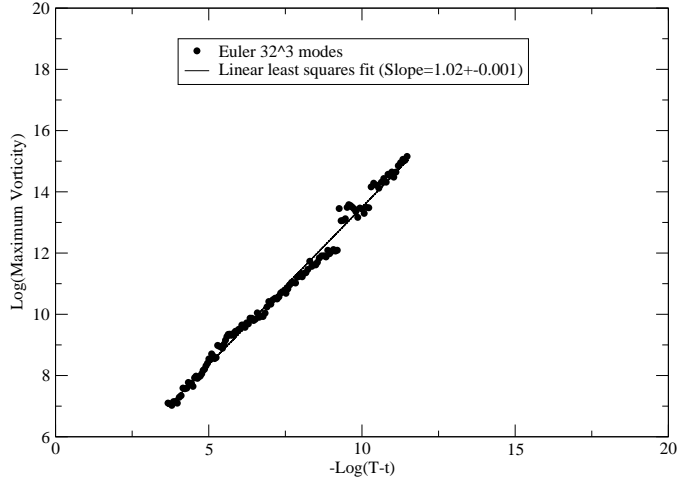}}
\caption{Log-log plot of the maximum vorticity versus $(T-t)^{-1}$ for the cycles 18-45 (about 10000 points)
  for the Euler equations with $32^3$ modes.} 
\label{fig_eul5}
\end{figure}

We also tried to see whether the solution near the purported blow-up satisfied
the Beale-Kato-Majda criterion for blow-up \cite{bkm}. The criterion states that if the maximal 
time of existence of a local solution is $T,$ then 
\begin{equation}
\label{bkm}
\int_{0}^{T} \|\omega \|_{L^{\infty}}(t) dt=\infty.
\end{equation} 
We present results from calculations where we tracked the largest component of the vorticity. 
The verification of the criterion is not straightforward because of the growth of round-off
error when the vorticity increases. 
We computed an 
approximation to the integral in equation (\ref{bkm}) by translating the maximum vorticity and the stepsizes during 
all the cycles to the original scale and then using the trapezoidal rule. We obtained for the integral in 
(\ref{bkm}) the value 464.38, which is not decisive for the finite time blow-up of 
$\int_{0}^{T} \|\omega \|_{L^{\infty}}(t) dt.$ Another way to look for a possible finite time blow-up is to 
assume that the maximum vorticity is $\sim (T-t)^{-\zeta},$ then plot $\log  \|\omega \|_{L^{\infty}}(t)$ versus 
$-\log (T-t)$ and find the linear least squares fit of the plotted function. 
Figure \ref{fig_eul5} is such a plot for cycles 18-45 together with the linear least squares fit. The slope of the linear 
fit is $\zeta=1.02\pm 0.001$ for the last 10000 values of the maximum vorticity. This suggests that the maximum vorticity is indeed 
behaving as $\sim (T-t)^{-\zeta}$ near $T=6.05,$ consistent with a finite-time blow up. It is not safe to make any stronger claims. On the 
other hand, the maximum value of the velocity (not shown here) remains of $O(1)$ during all the cycles, as one should expect.

\begin{figure}
\centering
\subfigure[]{\epsfig{file=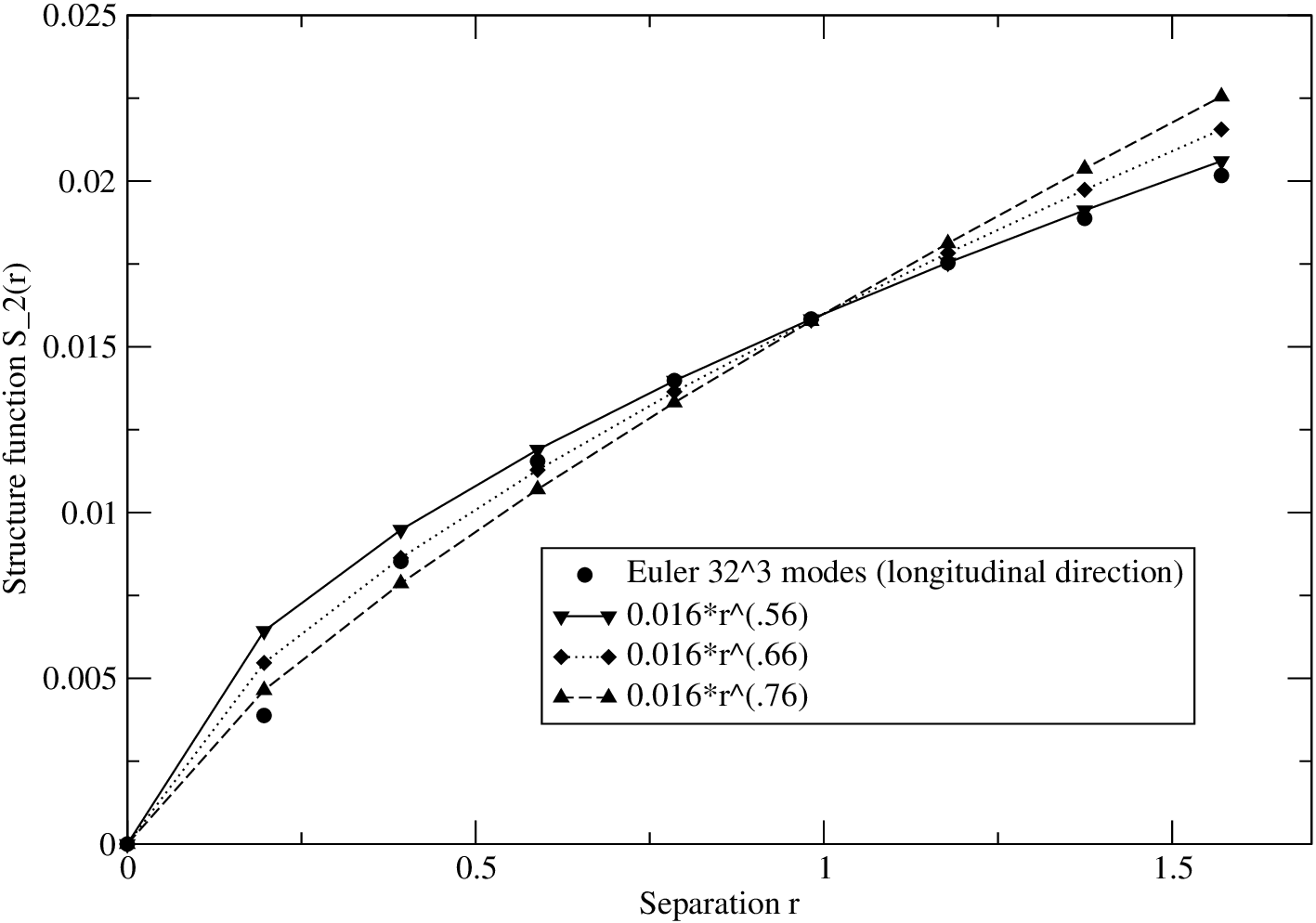,height=3.3in}}
\caption{Structure function of order 2 averaged over 45 cycles for the Euler equations
with $32^3$ modes. The power laws with the Kolmogorov exponents are inside an envelope of power laws with
lower and higher exponents.}
\label{fig_eul4}
\end{figure}

\begin{figure}
\centering
\subfigure[]{\epsfig{file=fig_eul14.eps,height=2.50in}}
\qquad
\subfigure[]{\epsfig{file=fig_eul15.eps,height=2.50in}}
\caption{Structure functions of order 3-4 averaged over 45 cycles for the Euler equations
with $32^3$ modes. The power laws with the Kolmogorov exponents are inside an envelope of power laws with
lower and higher exponents.}
\label{fig_eul41}
\end{figure}

\begin{figure}
\centering
\subfigure[]{\epsfig{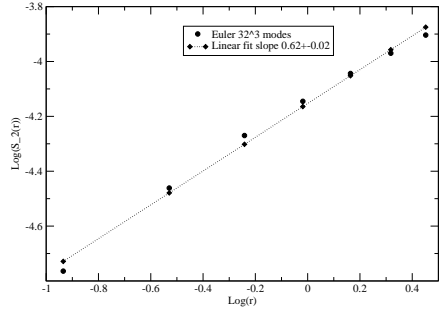}}
\qquad
\subfigure[]{\epsfig{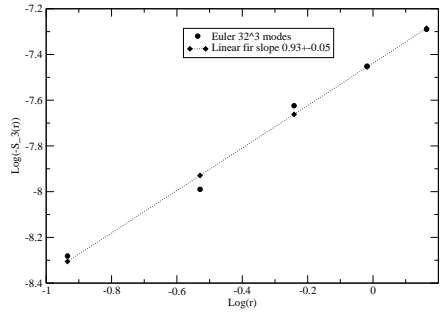}}
\qquad
\subfigure[]{\epsfig{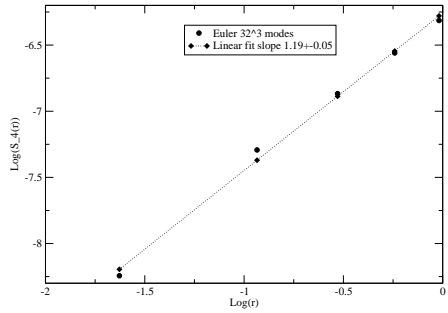}}
\caption{Log-log plot of the structure functions of order 2-4 averaged over 45 cycles for the Euler equations
with $32^3$ modes along with their least-sqaures fits.}
\label{fig_eul7}
\end{figure}

Figure \ref{fig_eul4} presents the averaged second order structure function for the Euler equations for the interval 
$[0,\frac{\pi}{2}],$ while Fig.\ref{fig_eul41}) shows the thrid and fourth order structure functions
with the same detail, but the results also corroborate our modest claims. 
We do not include the 0th cycle in the averaging of the structure functions, which can be 
thought of again as an equilibration step, helping to forget specifics of the initial condition. This omission 
is  consistent with our algorithm, which aims to reveal the generic structure of regions with 
highest vorticity and not problem dependent parameters. As in the case of Burgers, the averaged structure functions can only reveal 
the common features of the different cycles' small scales, since the large scales features differ from cycle to cycle. 
We have included in the figures power laws of the form $\alpha_p r ^{\beta_p},$ where $p$ is the order of the structure function. 
The figures include the power law predicted by Kolmogorov's theory, i.e. $\frac{p}{3}$ as well as 
the power laws with exponents $\frac{p}{3}\pm 0.1.$ 
We present this envelope of power laws to show that a calculation with $32^3$ modes does not allow an accurate
determination of the exponents. However, 
the range of values around the Kolmogorov values that are compatible with the numerical results is not broad. This
observation leads us to hope that the inertial range exponents for the Euler equation can be calculated
accurately when, in the future, we use our method with a larger number of Fourier modes (see also discussion in
Section \ref{conc}). For the sake of completeness, we include in figure \ref{fig_eul7} log-log plots and the corresponding least squares fits for the 
averaged structure functions of order 2-4. We see that the slopes of the fits are within the envelope of power laws shown in 
figures \ref{fig_eul4},\ref{fig_eul41}. In particular, we obtain the slope $0.62 \pm 0.02$ for the second order, $0.93 \pm 0.05$ for the 
third and $1.2 \pm 0.04$ for the fourth. The result for the fourth order structure function is marginally inside the envelope presented in 
figure \ref{fig_eul41}b, but this is to be expected. As one goes to higher order moments, the inadequacy of the resolution results in the 
envelope of power exponents to broaden. This fact prohibits us to make any stronger claims at this moment.

\begin{figure}
\centering
\subfigure[]{\epsfig{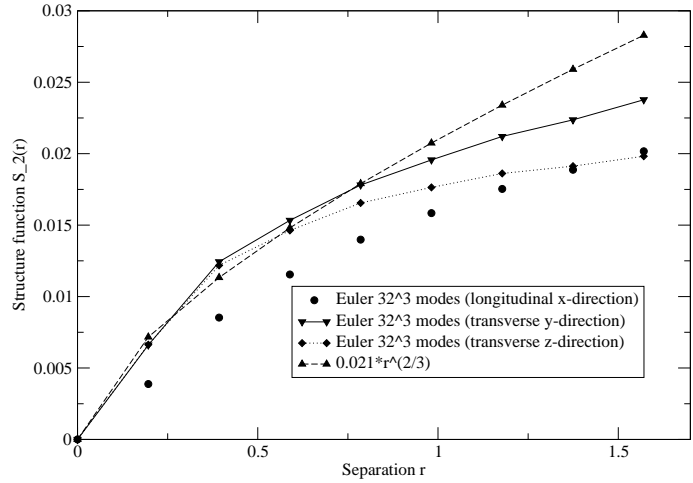}}
\caption{ Longitudinal and transverse second order structure functions along with Kolmogorov
 power law for transverse functions.}
\label{fig_eul42}
\end{figure}

\begin{figure}
\centering
\subfigure[]{\epsfig{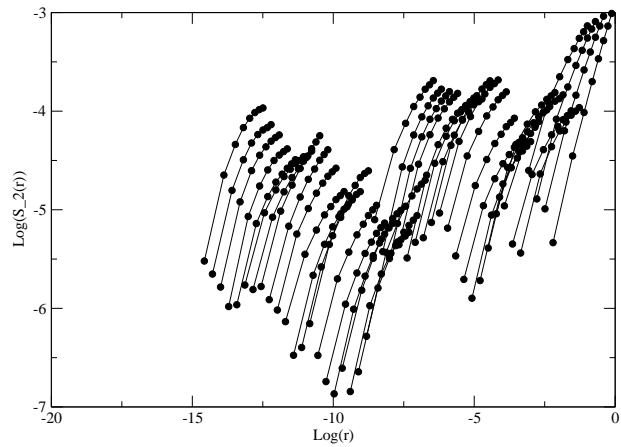}}
\caption{ Log-log plot of the second order structure function for the different cycles.
The structure functions for the different cycles are translated to the original scale.}
\label{fig_eul6}
\end{figure}

Figure \ref{fig_eul42} shows the longitudinal and transverse second order structure functions. If the flow is 
isotropic and a power law behavior holds for the second order structure function, then the longitudinal and transverse second structure 
functions should exhibit the same exponent, but with a different prefactor \cite{batchelor}. 
Experiments for high Reynolds numbers (e.g. \cite{sreeniv}) 
show that there is some discrepancy between the 
exponents of the longitudinal and transverse structure functions. This discrepancy is usually attributed to the the fact that 
perfect isotropy cannot be obtained in an experiment. In our numerical experiments we see that the longitudinal and transverse second structure 
functions do exhibit the same scaling behavior for small $r$ (to within the accuracy afforded by the numerics).

Finally, figure \ref{fig_eul6} is a log-log plot of the second order structure function for different
cycles; as before, the structure functions for all the cycles was converted to the original spatial scale. Suppose that the 
structure functions for the different cycles (translated to the original scale) are $\phi_0(r),\phi_1(r),\phi_2(r),\ldots$ If the structure 
function exhibits power law behavior, then the structure functions for the different cycles should be given by 
$\phi(r),\phi(r/\frac{4}{3}), \phi(r/(\frac{4}{3})^2),\ldots,$ where $\phi(r)$ is the power law that holds across all scales. 
To check that, plot $\phi_0(r)$ in log-log coordinates for $r\in (0,\Lambda],$ where $\Lambda$ is some prescribed interval 
where the scaling law is expected to hold. Plot $\phi_1(r)$ in log-log coordinates for $r\in (0,\Lambda/\frac{4}{3}],$ 
 $\phi_2(r)$ in log-log coordinates for $r\in (0,\Lambda/(\frac{4}{3})^2]$ etc. If the same scaling law holds across the cycles, the 
plots for the different cycles are parallel to one another, at a distance $\log(\frac{4}{3})$ apart.
 For the numerical experiments with the Euler equations, 
this is not exactly so. The reason is,
 that in three dimensions there is no longer a localized singular structure (like the shocks in one dimension). Thus, tracking of only 
the point of highest vorticity (and considering the inevitable inaccuracy of the numerical calculations) can result in jumping 
around different points of the singular structure. These different points can exhibit the same scaling behavior but with different 
numerical prefactors. This leads to the log-log plot of the second structure function being 
divided into clusters of cycles with different heights. However, the ranges of heights of the various clusters are 
in the same order of magnitude.


\section{Conclusions}\label{conc}

We have presented an algorithm that combines successive scaling with a spectral method in an attempt to 
probe the small scale structure of turbulence. 
We addressed the controversial issue of finite time blow-up for the solution of the Euler equations starting 
from the Taylor-Green vortex initial condition. 
We find that the behavior is consistent with a finite time blow-up 
of the vorticity while not being in any way desicive.

We used the algorithm to compute low order structure functions. For small distances, 
the structure functions exhibit power-law behavior. For the Euler equations, the Kolmogorov estimates of the power-law 
exponents are compatible with our results but we 
cannot conclude whether there exist corrections to the Kolmogorov estimates. Even though 
our calculations allow us to probe very fine scales, we need higher resolution, 
e.g. $512^3$ or $1024^3$ Fourier modes. Such resolutions  should be feasible through  parallelization and we 
expect to report on such calculations in the future. The merit of the algorithm is that
it reduces the resolution needed to decide the values of the exponents from something
astronomical to something merely very difficult. Before we attempt such calculations
we expect to improve our constructions, in particular in their treatment of boundary
conditions, as we have explained.

\section{Acknowledgements} We are grateful to Prof. Barenblatt, G.I., 
Prof. V.M. Prostokishin and Dr. Yelena Shvets for many helpful discussions and comments.
This work was supported in part by the National Science 
Foundation under Grant DMS 04-32710, and by the Director,
Office of Science, Computational and Technology Research, U.S.\ Department of 
Energy under Contract No.\ DE-AC03-76SF000098.


\begin{thebibliography}{99}

\bibitem{smith}
Smith, L.M. (1998), Ann. Rev. Fluid Mech. 30,
pp. 275-310.

\bibitem{chorin12}
Chorin, A.J. (1981),
Comm. Pure App. Math. 34  pp. 853-866.

\bibitem{barenblatt}
Barenblatt, G.I., Scaling. Cambridge University Press, Cambridge, 2002.


\bibitem{goldenfeld}
Goldenfeld, N. (1992),  Lectures on Phase Transitions and the Renormalization Group,
Perseus Books, Reading, Mass..


\bibitem{berger}
Berger, M. and Kohn, R. (1988),
Comm. Pure Appl. Math. 41, pp. 841-863.


\bibitem{landman}
Landman, M.J.,  Papanicolaou, G.C., Sulem, C.,  and Sulem, P. (1988),  Phys. Rev. A 38
pp. 3837-3847.


\bibitem{frederico}
Braga, G.,  Furtado, F. and Isaia, V. (2004),
Proc. Fifth Int.
Conf. Dynamical systems and Partial Differential Equations, Pomona CA, pp. 1-12.


\bibitem{gear}
Gear, C.W.,  Li, J.,  and Kevrekidis, I. (2003),
Phys. Lett. A 316 pp. 190-195.


\bibitem{kolm}
Kolmogorov, A.N. (1941),
Dok. Akad. Nauk USSR 30, pp.299-312.



\bibitem{batchelor}
Batchelor, G.K. (1960), The Theory of Homogeneous Turbulence, Cambridge University Press.


\bibitem{bernard}
Bernard, P. and Wallace, J. (2002), Turbulent Flow: Analysis, Measurement and Prediction,
Wiley, Hoboken NJ.

\bibitem{monin}
Monin, A., and Yaglom, A. (1971), Statistical Fluid Mechanics, MIT press, Cambridge, MA.




\bibitem{bcbc}
Barenblatt, G.I. and Chorin, A.J. 
(1998),
Proc. Symposia Appl. Math. AMS 54 ,
pp. 1-25.

\bibitem{BC04}
Barenblatt, G.I. and Chorin, A.J. (2004), , 
Proc. Nat. Acad. Sc. USA 101, pp. 15023-15026.


\bibitem{kadanoff}
Kadanoff, L. (2002),  Statistical Physics: Statics, Dynamics, and Renormalizatiion,
World Scientific, Singapore.

\bibitem{CH}
Chorin, A.J. and Hald, O. (2005), Stochastic Tools for Mathematics and Science,
Springer, NY.


\bibitem{spalart}
Spalart, P.R. (1988),
Fluid Dynamics of Three-Dimensional Turbulent Shear Flows and Transition,
AGARD-CP-438, AGARD, Neuilly-sur-Seine, France, pp. 5.1-5.13.



\bibitem{nordstrom}
Nordstr\"{o}m, J., Nordin, N.,  and Henningson, D.S. (1999),
SIAM J. Sci. Comput. 20 pp. 1365-1393.



\bibitem{hair}
Hairer, E.,  N\"orsett, S.E.,  and Wanner, G. (1987),
Solving Ordinary Differential Equations I-II, Springer, NY.

\bibitem{canuto}
Canuto, C., Hussaini, M.Y.,  Quarteroni, A.  and  Zang, T.A. (1988),
Spectral Methods in Fluid Dynamics, Springer, NY.


\bibitem{chorin0}
Chorin, A.J. (1969),
Math. Comp. 23 pp. 341-353.


\bibitem{chorin4}
Chorin, A.J. (1994), Vorticity and Turbulence, Springer, NY.



\bibitem{bkm}
Beale, J.T., Kato, T. and Majda, A., (1984),
Comm. Math. Phys. 94 , pp. 61-66.

\bibitem{morf}
Morf, R.H., Orszag, S.A. and Frisch, U., (1980), Phys. Rev. Lett. 44, pp. 572-575.


\bibitem{sreeniv}
Dhruva, B., Tsuji, Y.  and Sreenivasan, K.R. (1997),
Phys. Rev. E 56, pp.R4928-R4930.



\bibitem{xconst}
Constantin, P. (1994), SIAM Review 36 
pp. 73-98.







\end{thebibliography}
\end{document}